\documentclass[12pt]{article}
 \usepackage{amsmath}
 \usepackage{amssymb}

\begin{document}
\newtheorem{theorem}{Theorem}
\newtheorem{proposition}[theorem]{Proposition}
\newtheorem{lemma}[theorem]{Lemma}
\newtheorem{corollary}[theorem]{Corollary}
\newtheorem{claim}[theorem]{Claim}
\newcommand{\weil}{Weil-Petersson }
\newcommand {\teich}{Teichm\"{u}ller }
\newcommand{\T}{{\mathcal T}}
\newcommand{\Tbar}{\overline{\mathcal T}}

\newenvironment{proof}{ {\sc {\bf Proof}}}{ {\sc {\bf q.e.d.}} \\}
\newenvironment{remark}{\noindent {\bf Remark}}{{\sc } \\}
\newenvironment{definition}{\noindent {\bf Definition}}{{\sc } }

\begin{center}
\LARGE{\weil geometry of Teichm\"{u}ller--Coxeter complex and its finite rank property}
\end{center}

\vspace*{0.05in}
\begin{center}

Sumio Yamada\\
\hfill \\
\end{center}

\begin{abstract}
Resolving the incompleteness of \weil metric on \teich spaces
by taking metric and geodesic completion results in two distinct
spaces, where the Hopf-Rinow theorem is no longer relevant due to the
singular behavior of the \weil metric.  We construct a geodesic completion
of the \teich space through the formalism of Coxeter complex
with the \teich space as its non-linear non-homogeneous fundamental domain.
We then show that the metric and geodesic completions both satisfy a
finite rank property, demonstrating a similarity with the non-compact
symmetric spaces of semi-simple Lie groups.
\end{abstract}

\section{Introduction}

Let $\Sigma$ be a closed topological surface of genus larger than one.  Throughout this paper
we will assume that $\Sigma$ is equipped with some hyperbolic metric.  In~\cite{Y2}, the author has studied
geometry of the Weil-Petersson completed \teich space $\Tbar$ of the surface $\Sigma$.  It was shown that the
metric completion $\Tbar$ is an NPC (non-positively-curved) space, also called a CAT(0) space, and the $\Tbar$ has
a stratification by its boundary sets.   Set-theoretically it is the so-called augmented \teich
space(~\cite{Ab}\cite{Be}.)  The singular behavior of \weil metric as the conformal structure degenerates was first
studied by H.Masur~\cite{Ma}, which led to the author's work~\cite{Y2}.

We recall the setting of the bordification $\Tbar$. We first let $\cal S$  be the classes of
homotopically nontrivial simple closed curves on the surface $\Sigma_0$
up to homotopy equivalence.  This set can be identified with the set of simple closed
geodesics on the surface with a hyperbolic metric. Then define complex of curves $C({\cal S})$
as follows.  The vertices/zero-simplices of $C({\cal S})$ are the elements of $\cal S$.
An edge/one-simplices of the complex consists of a pair of homotopy classes of disjoint simple
closed curves. A $k$-complex consists of $k+1$ homotopy classes of mutually disjoint
simple closed curves.  A maximal set of mutually disjoint simple closed curves, which produces
a pants decomposition of $\Sigma_0$, has $3g-3$ elements.  We say a simplex $\sigma$ in
$C({\cal S})$ {\it precedes} a simplex $\sigma'$ provided $\sigma \subset \sigma'$, and we write $\sigma \geq
 \sigma'$. We say a simplex $\sigma$ in  $C({\cal S})$ {\it strictly precedes} a simplex $\sigma'$ provided
$\sigma \subsetneq \sigma'$, and write $\sigma > \sigma'$. This defines a partial
ordering by reverse inclusion in the complex of curves $C({\cal P})$, and thus makes it a partially ordered set
(poset.)  We define the null set to be the $(-1) $-simplex.  Then there is a $C({\cal S}) \cup \emptyset$--valued
function $\Lambda$, called labeling, defined on $\Tbar$ as follows.
Recall a point $p$ in $\T$ represents a marked Riemann surface $(\Sigma, f)$ with an
orientation-preserving homeomorphism $f: \Sigma_0
\rightarrow \Sigma$.  The \weil completion $\Tbar$ consists of bordification points
of $\T$ so that $\Sigma$ is allowed to have nodes, which are geometrically
interpreted as simple closed geodesics of zero length. Thus a point $p$ in $\Tbar \backslash
\T$ represents a marked noded Riemann surface $(\Sigma, f)$ with $f: \Sigma_0
\rightarrow \Sigma$. We now define $\Lambda(p)$ to be the simplex of free homotopy classes on
$\Sigma_0$ mapped to the nodes on $\Sigma$. We denote the fiber of
$\Lambda: \Tbar \rightarrow C({\cal S}) \cup \emptyset$ at a point $\sigma \in C({\cal S})$ by ${\T}_{\sigma}$.
We denote its \weil completed space by ${\Tbar}_\sigma$.
The completed space $\Tbar$ has the stratification
\[
\Tbar = \cup_{\sigma \in C({\cal S})} {\T}_\sigma
\]
where the original \teich space $\T$ is expressed as ${\T}_\emptyset$.  For each
collection $\sigma \in C({\cal S})$, each boundary \teich space ${\T}_\sigma$ is a geodesically convex subset of
$\Tbar$~\cite{Y2} . Here geodesic convexity means that given a pair of points in ${\T}_\sigma$, there is a distance-realizing
\weil geodesic segment connecting them lying entirely in ${\T}_\sigma$.   When $C=\emptyset$
this result corresponds to the \weil geodesic convexity of \teich space $\T$ first obtained by~\cite{W1} where the
geodesic length functional was shown to be convex.  Later a large family of convex functionals were found in~\cite
{Y1}.

The \weil metric defined on $\T$ is a smooth Riemannian metric whose sectional curvature is negative everywhere.
Hence the only flats (i.e. isometric embedding of Euclidean space ${\bf R}^n, n \geq 1$) are the \weil geodesics.
The fact that there is no strictly negative upper bound for the sectional curvature is explained by the fact that the
\weil completion $\Tbar$ does have higher dimensional flats (see the concluding remark in~\cite{Y2}, as well as
Proposition 16 in~\cite{W2}.)   Those flats arise when the collection $\sigma$ of mutually disjoint simple closed
geodesics separates the surface $\Sigma$ into multiple components.  Then the frontier set ${\T}_\sigma$
is a product space of the \teich spaces of the components separated by the nodes.  The number of the connected components is bounded by
$g+(\big[\frac{g}{2}\big]-1)$, which is achieved when the surface $\Sigma_\sigma$ is a union of $g$ once punctured tori,
and $\big[\frac{g}{2}\big]-1$ four-times punctured spheres. One can construct isometric embedding of Euclidean spaces
of dimension up to $g+(\big[\frac{g}{2}\big]-1)$ by considering a set of \weil geodesic lines, whose existence is
established in~\cite{DW}\cite{W2}, in different components of
the product space ${\T}_\sigma$.  In this sense the \weil completion $\Tbar$ is a space of {\it finite rank} where the
rank is bounded by $g+(\big[\frac{g}{2}\big]-1)$.

In this paper, we will use another definition of rank, which we call  FR, as first appeared in~\cite{KS2}.
\\

\begin{definition}\label{FR} An NPC (${\rm CAT}(0)$) space $(X, d)$ is said to be an {\bf FR} space if there exist $\varepsilon_0 >0$ and
$D_0$ such that any subset of $X$ with diameter $D > D_0$ is contained in a ball of radius $(1-\varepsilon_0)D/\sqrt{2}$.
\end{definition}

We make several remarks about FR spaces.  The definition can be interpreted as follows. If $X$ is an FR space, then among
all the closed bounded convex sets $F$ in $X$ with its diameter larger than $D_0$, there exists some positive integer $k$
such that
\[
\inf_{F \subset X} \frac{D(F)}{R(F)}\geq \sqrt{2} \sqrt{\frac{k+1}{k}}
\]
where $D(F)$ and $R(F)$ are the diameter and the circumradius of $F$ respectively.
It is well known that ${\bf R}^k$ is FR with the optimal/largest choice of $\varepsilon_0
= 1- \sqrt{k/k+1} >0$ which is realized by the standard k-simplex.  An infinite dimensional
Hilbert space is not an FR space, while a tree is an FR space with $\varepsilon_0
= 1- \sqrt{1/2} >0$.  It was shown in~\cite{KS2} an Euclidean building is an FR space with
$\varepsilon_0 = 1- \sqrt{k/k+1} >0$ with $k$ the dimension of chambers.  A CAT(-1) space
(e.g. hyperbolic plane ${\bf H}^2$) is FR with $\varepsilon_0$ which can be made
arbitrarily close to $1-
\sqrt{1/2}$   by taking the value of $D_0$ large.  Heuristically the number $\varepsilon_0>0$ detects the maximal
dimension of flats inside the space $X$, that is, the rank of the given NPC space.

We show that

\begin{theorem}
The \weil completion $\Tbar$ of a \teich space $\T$ is FR.
\end{theorem}

The proof of the theorem determines a lower bound of $\varepsilon_0$ to be
$1-\sqrt{k/k+1}$ with $k =6g-6$, which is larger (for $g>1$) than the {\it maximal} dimension
of the flats as described above, which was $g+(\big[\frac{g}{2}\big]-1)$.  The particular value of
$k$ here should be regarded as the maximal dimension of flat.  Namely if
there is a flat, its dimension cannot exceed $6g-6$.    On the other hand, the existence of those flats arising from
the product structure of the frontier sets of $\Tbar$ does not necessarily imply that the space is FR, as there may be
infinite dimensional flats
elsewhere.  The definition of the FR spaces utilizes only the convex
property of the distance function to describe the finite dimensionality of possibly existent flats,
without directly dealing with the singular behavior of the \weil metric tensor near the frontier set $\partial {\Tbar}$.
The statement of the theorem says that despite the lack of local compactness near the frontier set $\partial {\Tbar}$,
the \weil completion $\Tbar$ of the \teich space $\T$ exhibits a finite rank characteristics.

The finite rank theorem is proven by considering an auxiliary space. We first show that this space is FR.
The space contains the \weil completion ${\Tbar}$ as a complete convex subset, and the
FR condition for the  auxiliary space is inherited to its subset $\Tbar$, providing the proof
of the theorem.

This auxiliary space by itself is of  independent interest.
We will leave the detailed construction in the following sections, but
present the basic idea here.  The theory of Coxeter group was developed in an attempt to understand
combinatorial and geometric characterizations of tilings of standard spaces such as ${\bf R}^n$,
$S^n$ and ${\bf H}^n$ by reflecting convex polygons across the sides of the polygons.  Each reflection
is an involution, and the group generated by all the reflections is called the Coxeter group.
Naturally the geometry of the vertices of the polygon comes into the structure of the Coxeter group.
Around 1960 Jacques Tits has introduced the notion of an abstract reflection group, which he
named ``Coxeter group" $(W, S)$, which is a group $W$ generated by a set of involutions $S$, and a collection
of relations among the involutions $\{ (ss')^{m(s, s')} \}$, where $m(s, s')$ denotes the order of
$ss'$ and the relations range over all unordered pairs $s, s' \in S$ with $m(s, s') \neq \infty$.
The data $[m(s, s')]_{(s, s')}$ can be regarded as a matrix, and is said to constitute a Coxeter matrix.
The linchpin connecting the \weil geometry of $\Tbar$ and the Coxeter theory is the following theorem of
Wolpert's. \\

\noindent {\bf Theorem}~\cite{W3} {\it Given a point $p$ in ${\T}_C \subset {\Tbar}$, which represent a nodal surface
$\Sigma_\sigma$ the Alexandrov tangent cone with respect to the \weil distance function is isometric to
${\bf R}^{|\sigma|}_{\geq 0} \times T_p {\T}_\sigma$, where ${\bf R}^{|\sigma|}_{\geq 0}$ is identified with the 
orthant in ${\bf R}^{|\sigma|}$ with the standard metric.} \\

The significance of this result in our context is that it describes the geometry of the vertices, when $\Tbar$ is seen as a
convex polygon, and specifies a particular choice of the Coxeter matrix.
Namely for each $\sigma$ with $|\sigma|=1$, one can {\it reflect}
$\Tbar$ across the totally geodesic stratum ${\Tbar}_\sigma$.  Now for $\tau = \sigma \cup \sigma'$
with  $\sigma$ and $\sigma'$ representing disjoint simple closed geodesics, the relation
$m(s_\sigma, s_{\sigma'}) =2$ has a geometric representation where four copies of $\Tbar$
can be glued together around a point $q \in {\T}_\tau$ to form a  space whose tangent cone
at $q$ is a union of four copies of  ${\bf R}^{2}_{\geq 0} \times T_p {\T}_\tau$ (each ${\bf R}_{\geq 0}^2$
is regarded as a quadrant in $(x,y)$-plane,) isometric to ${\bf R}^2  \times T_p {\T}_\tau$ on which the reflections $s_\sigma, s_{\sigma'}$ act as reversing of the orientations of the $x, y$ axes for ${\bf R}^2$.

Hence the Coxeter group $(W, S)$, where the generating set $S$ is identified with the elements of
${\cal S}$, and the relations are specified by the Coxeter matrix  whose components satisfying
i) $m_{ss}=1$, ii) if $s \neq s'$, and if there is some simplex $\sigma$ in $C({\cal S})$
containing $s$ and $s'$, then define $m_{ss'}=2$, and
iii) if $s \neq s'$, and if the geodesics representing $s$ and $t$ intersect on $\Sigma_0$ then
$m_{ss'} = \infty$,
has a {\it geometric realization} acting on a collection of copies of $\Tbar$'s.  The Coxeter group
with such a Coxeter matrix as above is said to form a {\it cubical} complex, for it can be realized geometrically
as each generating element can be represented as a linear orthogonal reflection across the face of a unit cube
in ${\bf R}^{|S|}$.   The remarkable phenomena here is that despite of the fact that the generating set is
infinite, we have a geometric realization of the Coxeter group $(W, S)$ action (which is very far from linear)
on a space modeled on a finite dimensional space $\T$, albeit the bordification $\Tbar$ encodes non-locally compact geometry due to the singular behavior of the \weil metric tensor.  Note that the singularity is
also manifest in the fact that the reflecting wall $\Tbar_\sigma$ with $|\sigma|=1$ is of complex codimension
one, instead of its being a real hypersurface as in the standard Coxeter theory.  We also make a remark that
\weil metric defined on each stratum ${\T}_\sigma$ is K\"{a}hler, hence the $D({\Tbar}, \iota)$ is
a ``simplicial" complex with each face equipped with a K\"{a}hler metric, a situation unattainable in
a  {\it real} simplicial complex, where the reflecting walls are real hypersurfaces.
The space obtained by the action of the Coxeter group on $\Tbar$,
which we will call development $D({\Tbar}, \iota)$, is then shown
to be ${\rm CAT}(0)$ via the Cartan-Hadamard theorem, and also to be
geodesically complete. 

Given a Riemannian manifold $M$, and a codimension-{\it two} submanifold $S \subset M$, the open manifold
$N:=M\backslash S$ has $M$ as its metric completion as well as the
geodesic completion with respect to the Riemannian distance
function. (Consider $M={\bf R}^2$, $S=\{0\}$ and $N$ the punctured
plane, for example.) Here the analogous picture is $N = \T$, $S=
{\Tbar} \backslash {\T}$, and $M$ is either $\Tbar$ or $D({\Tbar},
\iota)$, depending on whether the completion is taken to be metric
or geodesic. The disparity, that $\Tbar$ is metrically complete
but not geodesically complete, is caused by the singular behavior
of the \weil metric near the strata, where the points in the
frontier set can be modelled as vertex points of cusps~\cite{W2}.

A consequence of the geodesic extention property of the development $D({\Tbar}, \iota)$ is that
the inverse map ${\rm exp}_p^{-1}$ of the ``exponential map" from $D({\Tbar}, \iota)$ to the
tangent cone isometric to ${\bf R}^n$ is  surjective,
as every geodesic segment starting at $p$ can be extended to a geodesic line so that the image by the inverse exponential map is  an entire real line through the origin of
${\bf R}^n$.  It is precisely this point that will be needed in the proof of the finite rank theorem.  Namely the
finite rank of the space $D({\Tbar}, \iota)$ is shown by using Caratheodory's theorem about convex sets in
${\bf R}^n$, which then implies the inequality between the diameter $D$ and the circumradius $R$ of convex sets
in $D({\Tbar}, \iota)$. Without the surjectivity of  the inverse exponential map  ${\rm exp}_p^{-1}$, the
correspondence between the space of geodesics and the space of directions breaks down.  Also one notes that the
geodesic completeness
is understood in the sense that any geodesic segment can be extended to {\it a} geodesic line, and there may be
more than one extension (in fact uncountably many extensions) making ${\rm exp}_p$ multi-valued.
This is once again due to the singular behavior of the \weil metric tensor near the frontier sets, where the sectional
curvature can blow down to $-\infty$, which causes that the behavior
of geodesics resembles that of geodesics in ${\bf R}$-trees.

Lastly we consider the existence question of harmonic maps whose target spaces are the development $D({\Tbar}, \iota)$
and  $\Tbar$.  There are two issues  concerning the existence question caused by the geometry of the spaces.
One is the fact that the spaces are not a CAT($\kappa$)
space for any $\kappa <0$.  The strict negativity of the sectional curvature $\kappa$ of the comparison model space is
known to help establishing existence (see Theorem 2.3.1 in~\cite{KS1}.)  The other is the fact as seen above that there is a
$\bf Z$-action generated by Dehn twists $\gamma_\sigma$ around a simple closed geodesic $\sigma$ near the frontier set ${\T}_\sigma$, which
causes the lack of local compactness (see the argument in~\cite{W2}, as well as  Theorem 2.2.1 in~\cite{KS1}.)   Note that
the FR condition requires the space to be neither CAT($\kappa$) with $\kappa<0$, nor locally compact.
Knowing the space to be FR rectifies the situation, and one has the following existence criteria~\cite{KS2}.
We let $X$ be either $D({\Tbar}, \iota)$ or ${\Tbar}$.

\begin{theorem} Let $\Gamma$ be the fundamental group of a compact Riemannian manifold $M$, and let $\rho$ be an isometric action of
$\Gamma$ on $X$.  Either there exists an equivalence class of rays, fixed by the
$\rho$-action, or there exists a $\rho$-equivariant harmonic
map $u:\tilde{M} \rightarrow X$, where $\tilde{M}$ is the universal covering space of $M$.
\end{theorem}

We remark here that this dichotomy is the exact
analogue of the much investigated situation
for harmonic map into non-compact symmetric spaces $G/K$ for semi-simple Lie group $G$.
The isometric action is induced by a representation of $\Gamma$ in the isometry groups of the
spaces ${\Tbar}$ and $D({\Tbar}, \iota)$.  The former group is the extended mapping class group
as shown by Masur and Wolf~\cite{MW}.  The latter group contains a group which is the semi-direct
product of the extended mapping class group and the Coxeter group, in which the extended mapping class group
is a normal subgroup.

The paper is organized as follows.  In the section {\bf 2}, we construct the Coxeter complex $D({\Tbar}, \iota)$
from the Coxeter matrix as described above.  Next we look into the local properties of $D({\Tbar}, \iota)$,
in particular note that the reflection is regarded as gluing of {\rm CAT}(0) spaces, and that the
basic construction of the Coxeter complex induces a nice \weil tangent cone structures.  In the section {\bf 4},
we combine the global and local properties of $D({\Tbar}, \iota)$ to show that the space is itself a {\rm CAT}(0)
space, which is also geodesically complete.  In the following section, we quote a theorem of Wolpert's~\cite{W2}
where a version of sequential compactness in the space of geodesics in $\Tbar$ is presented, using
concatenations of geodesic segments.  We see that  how those concatenations correspond to
geodesics in the development $D({\Tbar}, \iota)$ in a canonical way.  In section {\bf 6}, we prove the finite rank
theorem using the convex analysis and the geodesic extension property of the space $D({\Tbar}, \iota)$.
Finally in the last section we introduce an isometry group acting on $D({\Tbar}, \iota)$, constructed as
a semi-direct product of the extended mapping class group and the Coxeter group.  In this semi-direct product
group, the Coxeter group is regarded as a ``Weyl group", permutating the reflecting walls, which
is labeled by the set ${\cal S}$.

The author would like to thank Scott Wolpert for helpful discussions.

\section{Global Construction}

The \weil metric completion of \teich space is not geodesically complete. Our current goal is to
construct a space which function as a geodesic completion of $\T$ with respect to the \weil
metric.

The setting we have is well-suited for the language of {\it stratified
space}~\cite{BH}: \\

\begin{definition}
$(\Tbar, \{{\Tbar}_{\sigma} \}_{\sigma \in C({\cal S})})$ is called
stratified set satisfying the following properties.  \\
1) $\Tbar$ is a union of strata ${\Tbar}_{\sigma}$,\\
2) if ${\Tbar}_{\sigma} = {\Tbar}_{\tau}$ then $\sigma= \tau.$ \\
3) if the intersection ${\Tbar}_{\sigma} \cap {\Tbar}_{\tau}$ of two strata
is non-empty, then it is a union of strata,\\
4) for each $p \in {\Tbar}$ there is a unique $\sigma (p) \in C({\cal S})$ such that
the intersection of the strata containing $p$ is ${\Tbar}_{\sigma(p)}$. \\
\end{definition}
  In this way, the space $\Tbar$ is a stratified space with strata $\{{\Tbar}_{\sigma}\}$
indexed by the elements of the poset $C({\cal P})$.

In this particular space, each stratum is a topological space and ${\Tbar}_\sigma \cap
{\Tbar}_\tau$ for each $\sigma$ and $\tau$ is a close subset of both ${\Tbar}_\sigma$ and
${\Tbar}_\tau$, and the topologies which ${\Tbar}_\sigma$ and ${\Tbar}_\tau$ induce
on ${\Tbar}_\sigma \cap {\Tbar}_\tau$ is the same.  We can then define a topology, called
weak topology induced from the metric topology with respect to the \weil distance function,
namely a subset of $\Tbar$ is closed if and only if its intersection with each stratum is
closed.  Note that for a given point there is a number $r>0$ such that $B_r(p)$
intersects only finite number of stratum, that is, the weak topology is locally finite.
Now for a given point $p$ in $\Tbar$ consider the set
\[
C_x = \cup_{s \notin {\cal S}(x) } \Tbar_s
\]
where $S(x) := \{ s \in {\cal S} | x \in {\Tbar}_s \}.$
As $\{{\Tbar}_s\}_{s \in {\cal S}}$ is a locally finite family of closed subspaces of $\Tbar$,
each $C_x$ is closed.  We call the complement of $C_x$ {\it open star} of $x$, and
denote it by ${\rm st}(x)$.

Next we demonstrate that one can associate a simple complex of groups~\cite{BH} to the stratified space $({\Tbar} , \{{\Tbar}_{\sigma} \}_{\sigma \in C({\cal S})})$.
We first define a Coxeter group $W$ with its generating set identified with ${\cal S}$.  Its defining
relations are of the form $(ss')^{m_{ss'}}$, one for each pair $s, s' \in {\cal S}$, satisfying\\
1) $m_{ss}=1$ and \\
2) if $s \neq s'$, and if there is some simplex $\sigma$ in $C({\cal S})$
containing $s$ and $s'$, then define $m_{ss'}=2$, \\
3) if $s \neq s'$, and if the geodesics representing $s$ and $t$ intersect on $\Sigma_0$
(namely there is no $\sigma$ in $C({\cal S})$ containing both $s$ and $s'$) then
$m_{ss'} = \infty$.   \\
The pair $(W, {\cal S})$ is called Coxeter system, and $M = [m_{ss'}]$
Coxeter matrix where $m_{ss'}$ is the $(s, s')$--component.

Given a subset $\hat{\cal S}$ of ${\cal S}$, we define $W_{\hat{\cal S}}$ be the Coxeter group with
generating set $\hat{\cal S}$ and relations $(uu')^{m_{uu'}}$ for $u, u' \in \hat{\cal S}$ where
$m_{uu'}$ is a component of the Coxeter matrix $M$ for $(W, \hat{\cal S})$. Then it is known that the natural
homomorphism $W_{\hat{\cal S}} \rightarrow W$ induced by the inclusion $\hat{\cal S} \hookrightarrow
{\cal S}$ is injective~\cite{Bo}. This implies that each element $u$ of $\hat{\cal S}$ generates
a cyclic group of order two, that if $u \neq u'$, then all the order-two cyclic subgroups are
distinct, and that for $u \neq u'$ with $m_{uu'} = \infty$, $u$ and $u'$ generate an infinite
dihedral group.

Recall that the complex of curves $C({\cal S})$ is a poset with a partial ordering.
It is a sub-poset of the poset consisting of all the subsets of $\cal S$ ordered by the reverse inclusion,
where the empty set $\emptyset$ is the maximal element in $C({\cal S})$.
For an element $\sigma$ of $C({\cal S})$,  we consider the subgroup $W_{\sigma}$ of $W$
generated by the elements $s$ of the vertex set $V(\sigma)$ with relations $(s s')^{m_{ss'}}$
for $s, s' \in V(\sigma)$. We denote the number of the elements in $V(\sigma)$ by $|\sigma|$.
By the way $m_{ss'}$ is chosen, the subgroup $W_{\sigma}$ is
isomorphic to the abelian group of order $2^{|\sigma|}$ acting on ${\bf R}^{|\sigma|}$ which permutes the
$2^{|\sigma|}$ orthants via reflections across the coordinate hyperplanes.

We introduce a concept which formulates a structure over a family of groups.\\

\begin{definition}
A simple complex of group $W(C({\cal S}))=(W_\sigma, \psi_{\tau \sigma})$
consists of the following;\\
1) for each $\sigma \in C({\cal S})$, a group $W_\sigma$, called the local group at $\sigma$
generated by vertex sets of $\sigma$ ;\\
2) for each $\tau < \sigma$, an injective homomorphism $\psi_{\tau \sigma} : W_{\sigma} \rightarrow
W_{\tau}$ such that if $\tau < \sigma < \rho$, then $\psi_{\tau \rho} = \psi_{\tau \sigma}
\psi_{\sigma \rho}.$  In particular for $\emptyset \in C({\cal S})$, $W_\emptyset$ is the
trivial group, and $W_{C({\cal S})} = W$.
\end{definition} \\

Associated to the simple complex of groups $W(C({\cal S}))$ over $C({\cal S})$, one has
the group
\[
\widehat{W(C({\cal S})) }:= {\lim_{\longrightarrow}}_{\sigma \in C({\cal S})} W_\sigma
\]
which is the direct limit, or equivalently the amalgamated sum of the system of groups and monomorphisms
$(W_\sigma, \psi_{\tau \sigma})$. It is obtained by taking the free product of the groups
$W_\sigma$ and making the identifications $\psi_{\tau \sigma} (h) = h$ for all $h \in W_\sigma$ and for any $\tau < \sigma$: namely for $W_\sigma= \langle {\cal A}_\sigma \, | \, {\cal R}_\sigma \rangle$ with ${\cal A}_\sigma$ the generators
and ${\cal R}_\sigma$ their relations, $\widehat{W(C({\cal S})}$ is presented as
\[
\langle \amalg_\sigma {\cal A}_\sigma \, | \, {\cal R}_\sigma, \, \psi_{\tau \sigma} (a) = (a), \forall a \in {\cal A}_\sigma,
\forall \tau < \sigma \rangle
\] 
The natural homomorphisms $\iota_\sigma :
W_\sigma \rightarrow \widehat{W(C({\cal S}))}$ give a canonical simple morphism
$\iota : W(C({\cal S})) \rightarrow \widehat{W(C({\cal S}))}$, and in general $\iota_\sigma$ is
not injective. In our case, however, the injective homomorphism  $\psi_{\tau \sigma}:W_\sigma \rightarrow W_\tau$ is the inclusion induced by the inclusion $\sigma \subset \tau$ of the generating sets.  Also note that   the amalgamated sum $\widehat{W(C({\cal S}))}$ is canonically isomorphic to the full Coxeter group $W$, as $ \widehat{W(C({\cal S}))} \subset W$ is clear, and the other inclusion follows by noting that all the relations among the generating sets ${\cal S}$ defining the Coxeter group $W$ are present in the set of the relations ${\cal R}_\sigma$ for $\sigma$ in $C({\cal S})$ defining $ \widehat{W(C({\cal S}))}$.   The result quoted above~\cite{Bo} guarantees that $\iota_\sigma$ is injective for each $\sigma$. 

Given  a stratified space $(X, \{X_\sigma\}_{\sigma \in {\cal P}})$ where
a group $G$ is acting by strata preserving morphisms, a subset
$Y \subset X$  is called a {\it strict fundamental domain} if it contains exactly one point from each orbit, and if
for each $\sigma \in {\cal P}$ there is a unique $p(\sigma) \in {\cal P}$ such that
$g.X_\sigma = X_{p(\sigma)} \subset Y$ for some $g \in G$.

Now we construct a space on which there is a group action with strict fundamental domain
$\Tbar$.     For this goal, we utilize so-called basic construction (\cite{BH} Theorem 12.18)
which is a way to construct a space $X$ from a fundamental domain $Y$ and
a group with a specified data about isotropy subgroups on each stratum.  Naively what we do
is to consider the product space $W \times {\Tbar}$, and to take a quotient space by equivalence relations
induced by isotropy groups.

Recall that we have the stratified space $({\Tbar}, \{{\Tbar}_{\sigma} \}_{\sigma \in C({\cal
S})})$ indexed by the poset $C({\cal S})$, the simple complex of groups
$W(C({\cal S}))= (W_\sigma, \psi_{\tau \sigma} )$ over $C({\cal S})$, and the canonical
simple morphism $\iota: W(C({\cal S})) \rightarrow W$, where $W$ is the Coxeter group
with the generating set identified with $\cal S$.
Let $D({\Tbar}, \iota)$ be the set which is the quotient of $W \times {\Tbar}$ by the
equivalence relation
\[
(g, y) \sim (g', y') \Longleftrightarrow y = y' \mbox{ and } g^{-1}g' \in W_{\sigma (y)}
\]
where ${\Tbar}_{\sigma (y)}$ denotes the smallest stratum containing $y$.  We write $[g,y]$
to denote the equivalence class of $(g, y)$.

Let a poset $D(C({\cal S}), \iota)$  be the
disjoint union of cosets $\amalg_{C({\cal S})} W \backslash W_{\sigma}$, whose element
is written as $(g W_\sigma, \sigma)$ with $\sigma \in C({\cal S})$ and $g W_\sigma \in
W \backslash W_\sigma$, a coset of $W_\sigma \subset W$. The partial ordering is given
by $(g W_\tau) < (g' W_\sigma)$ if and only if $\tau < \sigma$ and $g^{-1} g' \in
W_\tau$.  $W$ acts on $D(C({\cal S}), \iota)$ naturally by $g'(g W_\sigma, \sigma) = (g'g W_\sigma, \sigma)$.
Note that the map $\sigma \mapsto (W_\sigma, \sigma)$ identifies $C({\cal S})$ with
a subposet of $D(C({\cal S}), \iota)$.  Then $D({\Tbar}, \iota)$ is a stratified set over
the poset $D(C({\cal S}), \iota)$ where the stratum indexed by $(g W_\sigma, \sigma)$ is written as
$\cup_{y \in {\Tbar}_\sigma} [g, y]$.

The group $W$ acts on $D({\Tbar}, \iota)$ by strata preserving automorphism according to
the rule $g' [g, y] = [g'g, y]$.  By identifying $y$ in $\Tbar$ with $[1, y]$
in $D({\Tbar}, \iota)$, $\Tbar$ is regarded as a strict fundamental domain of the action.

Now we resort to a statement (\cite{BH} Proposition 12.20) to see the space $D({\Tbar},
 \iota)$ is simply connected. Namely the facts that ${\Tbar}$ is simply connected,
that each stratum of ${\Tbar}$ is arcwise connected, and that the canonical simple 
morphism is injective for each $\sigma$, implies that the
development $D({\Tbar}, \iota)$ is simply connected.  We remark that in the proof of this
result it is essential that the canonical morphism $\iota$ induces an injective homomorphism
$\iota_\sigma: W_\sigma \rightarrow W$ for each $\sigma \in C({\cal S})$.

\section{Local Construction}

Recall~\cite{BH} that an NPC (CAT(0)) space is a complete metric space $(X, d)$ which satisfies the following two hypotheses:

\begin{itemize}
 \item (Length Space) For any two points $x_0, x_1$ in $X$, there is a rectifiable curve $\gamma$ from $x_0$ to $x_1$
such that
\[
d(x_0, x_1) = {\rm Length}(\gamma).
\]
We call such a curve a geodesic.
 \item (Triangle Comparison) Given any three points $z, x_0, x_1$ in $X$, $\lambda$ in $[0,1]$ and a geodesic $\gamma$
from $x_0$ to $x_1$, let $x_\lambda$ denote the point which is a fraction $\lambda$ away from $x_0$ to $x_1$ along $\gamma$.
Write
\[
d(z, x_0) = d_0,  d(z, x_1)=d_1, d(z, x_\lambda) = d_\lambda, d(x_0, x_1)=L
\]
For an ${\bf R}^2$ comparison triangle of side lengths $d_0, d_1$ and $L$, we require that $d_\lambda$ less than or equal to
the distance from the vertex corresponding to $z$ and the point fraction $\lambda$ of the way from the point corresponding to
$x_0$ to the point corresponding to $x_1$.  The precise inequality is written as
\[
d_\lambda^2 \leq (1-\lambda) d_0^2 + \lambda d_1^2 - \lambda(1-\lambda) L^2
\]
\end{itemize}

We remark here that the (metric) completeness of the ${\rm CAT}(0)$ metric space is not always part of the definition, and at times it is more convenient
not requiring it.

We will demonstrate that the \weil completion $\Tbar$ can be extended to a
bigger ${\rm CAT}(0)$ space
by a gluing construction.   We quote the following
result which allows a general
gluing of two ${\rm CAT}(0)$ spaces along a complete convex subset.  \\

\noindent {\bf Theorem}(Y. Reshetnyak~\cite{Re}, Bridson-Haefliger~\cite{BH})  {\it Let $X_1$ and $X_2$ be ${\rm CAT}(0)$ spaces (not necessarily complete) and let $A$ be a complete metric space.  Suppose that for $j=1,2$, we are given isometries $i_j : A \rightarrow A_j$. where
$A_j \subset X_j$   is assumed to be convex.  Then  $X_1 \sqcup_A X_2$ is an ${\rm CAT}(0)$ space.} \\

Here the space $X_1 \sqcup_A X_2$ is the quotient space of the disjoint union $X_1 \coprod X_2$ by the equivalence relation
generated by $i_1(a) \sim i_2(a)$ for all $a \in A$.  The resulting space, called gluing space or amalgamation of the spaces $X_1$ and $X_2$ along $A$ has a canonical distance between $x \in X_j$ and $y \in X_{j'}$ defined as follows:
\begin{eqnarray*}
d(x, y) = &  d_j(x, y) &  \mbox{ if } j = j' \\
d(x, y) = &  \inf_{a \in A} \{ d_j (x, i_j(a)) + d_{j'} (x, i_{j'}(a)) \} &  \mbox{ if } j \neq j'
\end{eqnarray*}

Now given the \weil completion $\Tbar$ of the \teich space $\T$ of genus $g$ surface $\Sigma$,
we can glue a new space ${\Tbar}_\sigma$ to the original space $\Tbar$ along  a complete convex subset ${\Tbar}_{\tau}$
with $\sigma \subset \tau$ (or equivalently partially ordered as $\sigma > \tau$.)  The resulting space is an NPC space,
according to the gluing theorem above.

Using this theorem, we glue another copy of $\Tbar$ to the original $\Tbar$ along ${\Tbar}_\sigma$, where $\sigma$ is
a zero simplex in $C({\cal S})$ hence representing a simple closed geodesic $s(\sigma) \in {\cal S}$ in $\Sigma$.
Note that the gluing set ${\Tbar}_\sigma$ is a complete convex metric space as needed for the application of the gluing theorem.

Next we define the tangent cone of an ${\rm CAT}(0)$ space $X$. We define the space of directions $S_O X$ at $O$ to be the quotient space
of the space of geodesics starting at $O$ with the equivalent relation $\sigma_0 \sim \sigma_1$ if the angle between them
is zero.  The tangent cone $C_O X$ of $X$ at $O$ is defined to be the image of the map $\exp^{-1}: X \rightarrow C_O X$ which we
define by
\[
\exp_O^{-1}(p) = (d(O, p), \omega (p))
\]
where $\omega (p)$ is the direction of the geodesic segment connecting $O$ and $p$.  By
definition, the map preserves the distance along the radial geodesic rays originating at $O$
.  The notation $\exp^{-1}_O$ is of course
suggestive of the inverse map of the exponential map $\exp_O: T_O M \rightarrow M$ defined on each point of a Riemannian
manifold $M$.  Recall that the exponential map on a Riemannian manifold is a local diffeomorphism, and
hence the equivalence relation among geodesics is trivial.  However one cannot define exponential maps on
${\rm CAT}(0)$ spaces, as there may be more than one way to {\it exponentiate} a given initial direction.  Note that
even then the map $\exp^{-1}_O$ is still well-defined.

We now quote the following result obtained by S.Wolpert, which describes the \weil tangent cone of \weil completion of
the \teich space ${\Tbar}$.

\noindent {\bf Theorem}~\cite{W3} {\it Given a point $p$ in ${\T}_C \subset {\Tbar}$, which represent a nodal surface
$\Sigma_\sigma$ the Alexandrov tangent cone with respect to the \weil distance function is isometric to
${\bf R}^{|\sigma|}_{\geq 0} \times T_p {\T}_\sigma$, where ${\bf R}^{|\sigma|}_{\geq 0}$ is identified with the 
orthant in ${\bf R}^{|\sigma|}$ with the standard metric.}

 Hence for a point $q \in {\T}_\sigma$ with
 $|\sigma|=1$, the tangent cone is ${\bf R}_{\geq 0} \times {\bf R}^{6g-6-2}$ with $T_q {\Tbar}_C$ identified as
 ${\bf R}^{6g-6-2}$. In order to ensure that the intersecting set of the two copies of the tangent cone $C_q {\Tbar}$
 is exactly $T_q {\Tbar}_\sigma$, it is enough to know that each $\Tbar$ is geodesically convex in the bigger space
${\Tbar} \sqcup_{{\Tbar}_\sigma} {\Tbar}$ and that all the geodesics with their tangent vectors at $q$ in $S_q
[{\Tbar}_\sigma]$ are locally exponentiated within the frontier set ${\Tbar}_\sigma$, which is a consequence of the 
face that every open geodesic segement is entirely contained in a single stratum~\cite{DW}\cite{W2}.
Therefore the tangent cone $C_q [{\Tbar} \sqcup_{{\Tbar}_\sigma} {\Tbar}]$ is identified with 
$({\bf R}_{\geq 0} \times {\bf R}^{6g-6-2}) \sqcup_{{\bf R}^{6g-6-2}} ({\bf R}_{\geq 0} \times {\bf R}^{6g-6-2})$. 
Note that the metric across the gluing hyperplane
$\{ (0, y) | 0 \in {\bf R}_{\geq 0} \cap {\bf R}_{\geq 0} , \,\,\, y \in  {\bf R}^{6g-6-2} \}$
is yet to be  determined.

Consider the situation that the point $q$ lies in an open geodesic segment $\overline{p_0 p_1}
$ where $p_0$ lies in the original \teich space $\T$ and $p_1$
in the glued \teich space $\T$.  By the definition of the distance between $p_0$ and $p_1$,
the first variation of the sum of the distances
measured from $p_0$ and $p_1$ by varying $q$ within ${\Tbar}_\sigma$ vanishes.  Now let
$\sigma(t)$ be a geodesic in the frontier set ${\T}_\sigma$ through $q$.  Then by the angle
formula which equates the Alexandrov angle with the one-sided first derivative at $t=0$ of
the distance function as introduced above, we have
\[
0=\frac{d}{dt} \Big[ d(\sigma(t), p_0) + d(\sigma(t), p_1) \Big] \Bigg|_{t=0} = -\cos \angle_q
 (\sigma(t), p_0) -\cos \angle_q (\sigma(t), p_1).
\]
This implies that the projection of the tangent vectors of the geodesic rays $\overline{qp_0}$
 and $\overline{qp_1}$ onto the linear space ${\bf R}^{6g-6-2}$ are of the same length, of
the opposite direction.  This in turn shows that the perpendicular components to the frontier
set ${\Tbar}_\sigma$ of the two direction vectors have the same length, as those two unit vectors
are $\pi$ apart in the space of directions $S_q [{\Tbar} \sqcup_{{\Tbar}_\sigma} {\Tbar}]$.
Therefore the \weil distance function defined on
$[{\Tbar} \sqcup_{{\Tbar}_\sigma} {\Tbar}]$  is now linearized at $q$ so that the tangent cone, identified
above with $({\bf R}_{\geq 0} \times {\bf R}^{6g-6-2}) \sqcup_{{\bf R}^{6g-6-2}} ({\bf R}_{\geq 0} \times {\bf R}^{6g-6-2})$
is isometric to the the Euclidean space ${\bf R} \times {\bf R}^{6g-6-2}$.

Here we associate a {\it reflection} across the frontier set  ${\Tbar}_\sigma$, acting as a {\it wall},  to the gluing
space defined above. Namely for a given point $p$ in the {\it original} $\Tbar$ (or in the glued
 copy), let $s_\sigma(p)$ be $p$ in the {\it glued} copy of $\Tbar$ (or in the original copy
respectively.)  Here recall that each copy of $\Tbar$ is defined as a set of equivalence
class of pairs $(f, \Sigma)$ where $f: \Sigma_0 \rightarrow \Sigma$ is defined with respect
to a reference surface $\Sigma_0$. Thus the two points $p$ and $s_\sigma(p)$ are both marked by the
same reference surface $\Sigma_0$ and $f:\Sigma_0 \rightarrow \Sigma$ in the respective copy of $\Tbar$.
The map $s_\sigma$ is thus defined on
${\Tbar} \sqcup_{{\Tbar}_\sigma} {\Tbar}$ and it is one-to-one and onto, and its fixed point set
is precisely ${\Tbar}_\sigma$.  Also it is of order two, namely $s_\sigma(s_\sigma(p)) = p$ for all $p$.

Now for a point $q$ in ${\T}_\sigma \subset {\Tbar}$, as the tangent cone
$C_q [{\Tbar} \sqcup_{{\Tbar}_\sigma} {\Tbar}]$ is isometric to an Euclidean space, any geodesic $\sigma$
with $q = \sigma(0)$ can be extended beyond $q$ by the following steps.  First map the geodesic
$\gamma$ to the tangent cone $C_q [{\Tbar} \sqcup_{{\Tbar}_\sigma} {\Tbar}]$.  The
image of the map is a ray $ t \gamma'(0) \,\,\, t \geq 0$ emanating from the origin. Extend the ray across the origin to a line through the origin.  The added ray represents an equivalent class  of geodesics originating at $q$, which share the tangent vector $- \gamma'(0)$.  Pick a representative $\tilde{\gamma}$ of the equivalent class $[\tilde{\gamma}]$.  Namely $\tilde{\gamma}'(0) = -\gamma'(0)$.  Without loss of generality we assume that all the geodesics are arc-length parameterized here. Then the union of the two geodesic rays  $\gamma$ and $\tilde{\gamma}$ is locally length minimizing, as its first variations
of the length as $q$ varies in the frontier set ${\Tbar}_c$ vanishes.  Thus we have constructed a local extension of the geodesic $\gamma$ beyond its endpoint $q$. Notice that the extension of $\gamma$ is far from unique, as the equivalent class $[\tilde{\gamma}]$ contains uncountably many representatives.

For a point $q$ in a frontier space ${\Tbar}_{\tau(q)}$ where ${\Tbar}_{\tau(y)}$ has been defined to be the largest stratum containing the point $q$, we suppose $V(\tau) = \cup_{i=1}^{|\tau|} \sigma_i, \,\,\, |\sigma_i|=1$,
namely $\tau$ is an element in $C({\cal S})$ representing a union of $|\tau|$ mutually disjoint
simple closed geodesics $\{ \sigma_i\}$, we have the following gluing
construction, which ensures any geodesic terminating at $q$ is locally extended beyond $q$ in the bigger space. Let
$\sigma_1$ be an element of $\tau$, and glue a copy of $\Tbar$ along ${\Tbar}_{\sigma_1}$
to obtain a CAT(0) space ${\Tbar} \sqcup_{{\Tbar}_{\sigma_1}} {\Tbar}$.  We note that this space
has a frontier set ${\Tbar}_{\sigma_2} \sqcup_{{\Tbar}_{\sigma_1 \cup \sigma_2}} {\Tbar}_{\sigma_2}$ which is a
CAT(0) space of its own right with respect to the induced distance function from
${\Tbar} \sqcup_{{\Tbar}_{\sigma_1}} {\Tbar}$.  This frontier space is thus a complete convex subset in ${\Tbar}
\sqcup_{{\Tbar}_{\sigma_1}} {\Tbar}$. For $i=2$, we prepare an identical copy of the previously obtained space
${\Tbar} \sqcup_{{\Tbar}_{\sigma_1}} {\Tbar}$, and glue this copy to the original copy
along the complete convex subset ${\Tbar}_{\sigma_2} \sqcup_{{\Tbar}_{\sigma_1 \cup \sigma_2}} {\Tbar}_{\sigma_2}$
and equip the new space with the distance function given by the gluing theorem above.  The resulting space is once again CAT(0).  Recall~\cite{W3} the tangent cone $C_q {\Tbar}$ is
isometric to $({\bf R}_{\geq 0})^{|\tau|} \times T_q {\Tbar}_C \equiv ({\bf R}_{\geq 0})^{|\tau|} \times
{\bf R}^{6g-6-2|\tau|}$.  As we iterate the gluing process for $i = 1,2,\cdots$, we
denote the resulting space by $X_i(q)$.  For example $\Big({\Tbar} \sqcup_{{\Tbar}_{\sigma_1}}
{\Tbar} \Big) \sqcup_{{\Tbar}_{\sigma_1 \cup \sigma_2}} \Big({\Tbar}
\sqcup_{{\Tbar}_{\sigma_1}} {\Tbar} \Big)$ is written as $X_2(q)$.  Then
 the tangent cones at $q$ to the new spaces $X_i(q)$'s are isometric to the following spaces;
\[
C_q [X_1(q)] \equiv {\bf R} \times ({\bf R}_{\geq 0})^{|\tau|-1} \times
{\bf R}^{6g-6-2|\tau|}
\]
\[
C_q [X_2(q)] \equiv
{\bf R}^2 \times ({\bf R}_{\geq 0})^{|\tau|-2} \times {\bf R}^{6g-6-2|\tau|}
\]
and
\[
C_q [X_{|\tau|}(q)] =  {\bf R}^{|\tau|} \times {\bf R}^{6g-6-2|\tau|}.
\]
Hence after the gluing is repeated $|\tau|$ times, the tangent cone is isometric to an entire Euclidean space, any geodesic in $\Tbar$
terminating at $q$ can be now extended beyond $q$ as demonstrated above.

Hence for any point $q$, we have come up with a method to extend the original space $\Tbar$
via successive gluing so that in the resulting extended space any geodesic terminating at $q$ can be extended beyond $q$.

\section{Local-to-Global Geometry}

The next task is to relate the global construction to the local one.
The main tool to achieve the goal is the Cartan-Hadamard Theorem.  As we have
established in the preceding two sections, we have a {\it big} space which is obtained
by tiling the fundamental domain with the use of the Coxeter group action, and {\it small}
negatively curved spaces which is obtained by gluing construction along the frontier sets.
We will see that the small ones fit naturally within the big one, and that the
Cartan-Hadamard Theorem would specify the CAT(0) geometry of the big space.

We further show that the space displays the following property.

\begin{definition}
A geodesic metric space $X$ is said to have the {\it geodesic extension property} if
for every local geodesic $c:[a, b] \rightarrow X$ with $a \neq b$, there exists
$\varepsilon>0$ and a local geodesic $c':[a, b+\varepsilon] \rightarrow X$ such that
$c'|_{[a, b]} = c$
\end{definition}\\

We start with the local geometry of $D({\Tbar}, \iota)$.

\begin{proposition}
The development $D({\Tbar}, \iota)$ is of curvature $\kappa \leq 0$.
\end{proposition}

Here a space of curvature $\kappa \leq 0$ is defined to be a locally CAT(0)
space, i.e. for every $x$ in the space, there exists $r_x > 0$ such that
the closed ball $B_{r_x} (x)$ is a CAT(0) space.

\begin{proof}
Let $q$ be a point in the development of $D({\Tbar}, \iota)$ stratified over
$D(C({\cal S}), \iota)$.  Suppose $y$ is written as $[g, y]$ for $g \in W$ and
$y \in \Tbar$.  As by multiplying $g^{-1}$ to $q$ one can translate $q$ to
$g^{-1} [g, y] = [1, y]$,  without loss of generality we may assume that
$q$ is of the form $[1, y]$.  Considering the set $\{ [1, y] : y \in {\Tbar} \}$,
the point $q$ is identified with a point $y$ in ${\Tbar}_{\tau(y)}$ in the strict fundamental domain $\Tbar$ embedded in $D({\Tbar}, \iota)$.  Here $\tau(y)$ is
an element of the poset $C({\cal S})$ so that ${\Tbar}_{\tau(y)}$
is the smallest stratum containing $y$. Recall the local construction in the previous section, where we constructed a space $X_{|\tau(y)|}$ by gluing $2^{|\tau(y)|}$ copies of
$\Tbar$ around $q$.  We observed that the resulting space is CAT(0)
and any geodesic terminated at $y$ can be extended beyond $y$.

We claim that the space
$X_{|\tau(y)|}$ is embedded into the development
$D({\Tbar}, \iota)$.  In particular we show that $X_{|\tau(y)|}$ is the closed star
${\rm St}(q)$. Recall that the open star ${\rm st}(q)$ is defined to be the complement of the union
of strata which do not contain $x$, and the closed star, its closure.

In $X_{|\tau(q)|}(q)$, there are $2^{|\tau(q)|} $ copies of $\Tbar$ gathered around $q$.
The set ${\tau(q)}$ is a union of ordered set of mutually disjoint simple closed
geodesics $\cup_i \sigma_i$.  The gluing of two identical copies of $X_{k-1} (q)$
to obtain $X_k (q)$ $(k = 2,\cdots, |\tau(q)|)$ comes with the reflection map across
the frontier set along which the gluing has been performed.  The gluing set for
$X_{k}(q)$ is identified as the collection of $2^{k-1}$ copies of
${\Tbar}_\theta$ where $\theta$ is $\sigma_1 \cup \cdots \cup \sigma_{k}$. The union of those
copies of the frontier sets acts as
the reflecting {\it wall}, only it is of {\it complex} codimension one. (Note that
at the tangent cone level, the reflecting wall/hyperplane is of real codimension one.)  The
reflection is denoted by $r_{\sigma_{k}}:X_{k} (q) \rightarrow X_{k} (q)$.
As the ordering of the elements $\{ \sigma_i \}$ of $\tau$ is arbitrary, there is a set of $|\tau(q)|$
reflections defined on $X_{|\tau(q)|}(q)$, each of which reflects $\Tbar$ across
${\Tbar}_{\sigma_l}$ for $l = 1, \cdots , |\tau|$.  Those reflections acting on $X_{|\tau|}(q)$
generate a group of order $2^{|\tau(q)|}$.  This group is clearly isomorphic to the
Coxeter group $W_{\tau(q)}$ generated by the elements
$s_{\sigma_1}, \cdots, s_{\sigma_{|\tau|}}$ of $\cal S$
and the space $X_{|\tau(q)|}(q)$ is identified as the union of the orbits of points of $\Tbar$
by the special subgroup $W_{\tau(q)}$ of $W$.  Hence this space $X_{|\tau|}(q)$ is identified
as the closed star ${\rm St}(q)$ in $D({\Tbar}, \iota)$, as it is the union of all the strata
which contains the (zero) simplex $q$. Note here that
$q \in {\Tbar} \subset D({\Tbar}, \iota)$ is the maximal fixed point set of
the group action $W_{\tau}$, where recall that the inclusion
${\Tbar} \subset D({\Tbar}, \iota)$ is defined by the
embedding $y \rightarrow [1, y]$.

By considering this identification map ${\rm St}(y) \rightarrow X_{|\sigma(y)|}$, one
can define the pull-back distance function on ${\rm St}(y) \subset D({\Tbar}, \iota)$,
which in turn makes the development $D({\Tbar}, \iota)$ a metric space of curvature
$\kappa \leq 0$, locally isometric to $X_{|\sigma(y)|}$ for each $y$.

\end{proof}

\begin{theorem}
The development $D({\Tbar}, \iota)$ is a CAT(0) space, with every non-constant
geodesic can be extended to a geodesic line ${\bf R} \rightarrow D({\Tbar}, \iota)$
.
\end{theorem}

\begin{proof}
The statement of the Cartan-Hadamard theorem for metric spaces, shown by M.Gromov~\cite{Gr},
W.Ballmann~\cite{Bl}, S.Alexander-R.Bishop~\cite{AB}, is as follows.  Though the full statement as it appears
in~\cite{BH} is more general, we present a version we need.

\noindent {\bf Theorem} {\it Let $X$ be a
complete connected metric space.  If $X$ is of curvature $ \kappa \leq 0$, then
the universal covering space $\tilde{X}$ is a CAT(0) space.}

Having known that the development $D({\Tbar}, \iota)$ is of curvature $\kappa \leq 0$,
and that it is simply-connected as well as (metrically) complete, we can apply the
Cartan-Hadamard Theorem to conclude that the space $D({\Tbar}, \iota)$ is a CAT(0)
space.  As for the geodesic extension property, the local construction have made
the space $D({\Tbar}, \iota)$ satisfy the geodesic extension property.
It is known~\cite{BH}(Lemma 5.8) that a complete geodesic metric space, which is CAT(0),
has the geodesic extension property if and only if every local geodesic $\sigma: [a, b]
\rightarrow X$ can be extended to a geodesic line.

\end{proof}

 \section{Geodesic Extension Property }

We start this section with a statement (Proposition 23)  in Wolpert's
paper~\cite{W2}: a compactness result of the space of geodesics.  
It has a natural interpretation when considered in the
context of ${\Tbar} \subset D({\Tbar}, \iota)$.  Below, given an element $\sigma$ of the complex of curves $C({\cal 
S})$, we denote the stabilizer of the stratum ${\Tbar}_\sigma$
by ${\rm Mod}_\sigma$, which is a rank $|\sigma|$ Abelian group consisting of products of
Dehn twists about the mutually disjoint simple closed geodesics
represented by the zero-simplices of $\sigma$.  \\

\noindent {\bf Theorem}\cite{W2} {\it Consider a sequence of unit-speed geodesics $\{ \theta_n(s)\}$ with initial
points converging to $p_0$, lengths converging to a positive value $L'$ , and parameter intervals converging
to $[t', t'']$ with $t''-t'=L'$.  Then, there exists a subsequence, which we denote by $\{ \theta_n(s) \}$ again,
convergent in the following sense.

There exist a partition $t'=t_0 < t_1 < \cdots < t_k = t''$ of the interval $[t', t'']$, and a sequence of points $p_0, \cdots, p_k$
with each $p_j \in {\Tbar}_{\sigma(p_j)}$, so that  $d(p_j, p_{j+1})=t_{j+1}-t_j$ for $j=0,\cdots k-1$ and the concatenation
of geodesic segments $\overline{p_0 p_1} \cup \overline{p_1 p_2} \cup \cdots \cup \overline{p_{k-1} p_k}$ is the unique
length minimizing curve connecting $p_0$ to $p_k$ and intersecting in order the strata $ {\Tbar}_{\sigma(p_1)}, \cdots
{\Tbar}_{\sigma(p_{k-1})}$.  Furthermore, the following conditions are met among the
elements $\sigma (p_j)$'s of $C({\cal S})$. Denote $\sigma(p_{j}) \cap \sigma (p_{j+1})$ by
$\tau_j$.  Then {\bf 1)} $\tau_0 > \sigma(p_1)$ (namely as simplices $\tau_0 \subsetneq \sigma(p_1)$)
if $k>1$, {\bf 2)} $\tau_{k-1} > \sigma(p_{k-1})$ if $k>1$, {\bf 3)}
$\tau_j < \sigma (p_j)$ and $\tau_j < \sigma (p_{j+1})$ for $j =1,\cdots , k-2$.

Then the subsequence $\{ \theta_n \}$ comes with sequences of products of Dehn
twists $T_{(j, n)} \in {\rm Mod}_{\sigma(p_j) \backslash \tau_j}, \,\,\, j=
0,\cdots, k-1$ such that on the parameter interval $[t_j, t_{j+1}]$ the geodesic
segments $T_{(j, n)} \circ \cdots \circ T_{(0, n)} \theta_n$ converge to
$\overline{p_j p_{j+1}}$ as $n \rightarrow \infty$
in the sense of parameterized unit-speed curves in $\Tbar$.
Furthermore the distance between the parameterized unit-speed curves $\theta_n$
and $\overline{p_0 p_{(k, n)}}$ with $p_{(k,n)} := (T_{(k-1, n)} \circ \cdots
\circ T_{(0, n)})^{-1} p_k$, tends to zero for $n$ tending to infinity.
The sequences of transformations $\{ T_{(0,n)} \}_{n=1}^\infty$ is either trivial
or unbounded. For $k >1$, the sequences of transformations
$\{ T_{(j,n)} \}_{n=1}^\infty, \,\,\,j = 1, \cdots , k-1$ is unbounded.         }
\\

The limiting path can be represented in the development $D({\Tbar}, \iota)$ as
the length minimizing geodesic starting at $[1, p_0]:=q_0$ and terminating at
$[\Pi_{j=0}^{k-1}(\Pi_{s \in {V(\sigma(p_j))}} s), p_k]:=q_k$ passing through
$[\Pi_{j=0}^{l-1}(\Pi_{s \in {V(\sigma(p_j))}} s), p_l]:=q_{l}, \,\,\, l=1,\cdots, k-2$.
Here recall each element $s$ of the vertex set $V(\sigma(p_j))$ is identified as a
generator of the abelian subgroup $W_{\sigma(p_j)}$ of the Coxeter group $W$.  The
commutativity of the generators $s$'s makes the product $\Pi_{s \in {V(\sigma(p_j))}} s$
well defined.
Note that this path is locally length minimizing due to the fact that it is a collection of
geodesic segments. The path cuts across the strata with the incoming angle equal to the outgoing
angle with respect to each stratum, which in turn means the vanishing of the first variation
of length as described in {\bf Local Construction} section.  The CAT(0) condition of
$D({\Tbar}, \iota)$ then ensures that it is globally length
minimizing, realizing the distance between the two end points.

A remark should be made here that the paragraph above is
suggesting to interpret the reflection $s_\sigma \,\,\, \sigma \in
V(C({\cal S}))$ as $``\lim_{n \rightarrow \pm \infty} (T_\sigma)^n"$.
Namely when the space of Fenchel-Nielsen deformations around $\sigma$
($\cong {\bf R}$), which contains the Dehn twists ($\cong {\bf Z}
\subset {\bf R}$),  is compactified to $S^1={\bf R} \cup \{\infty\}$ 
via the Stone-Cech compactification, the reflection across the wall is 
regarded as the infinity element.

Let ${\cal A}_{p_0}(L)$ be the set whose elements are geodesics
in ${\Tbar}$, of a fixed length $L$ originating at a fixed
point $p_0$ in ${\Tbar}$ terminating at a point in $\Tbar$.
The compactness theorem of Wolpert's above together with the representation of
the limiting concatenated path as a geodesic segment in $D({\Tbar}, \iota)$  say
that given an infinite sequence $\{ \theta_n\}$ in
${\cal A}_{p_0}$, there exists a subsequence whose limit is represented as a
length minimizing path $\overline{[1,p_0] q_k}$ in $D({\Tbar}, \iota)$,
where $q_k$ is as defined above.

Conversely for each point $q$ in $D({\Tbar}, \iota)$, the geodesic segment
$\overline{[1,p_0] q}$ in $D({\Tbar}, \iota)$ can be identified as a limit of
a sequence of geodesics in $\Tbar$ given by a
concatenation $\overline{p_0 p_1} \cup \cdots \cup \overline{p_{k-1} p_k}$ and
unbounded sequences of products of Dehn twists $T_{(j, n)}$, where the points $\{ p_l \}$
are specified by the previously defined relation
$[ \Pi_{j=0}^{l-1}(\Pi_{s \in V(\sigma (p_j))} s), p_l] = q_l$ where $q_l$'s are the
points the geodesic segment $\overline{[1, p_0] q}$ is passing the stratum of $D({\Tbar},
\iota)$ through, and $T_{(j, n)}$ are specified to be products of Dehn twists in
${\rm Mod}_{\sigma(p_j) \backslash \tau_j}$, unbounded in $n$.

Now denote by $\tilde{{\cal A}}_{[1,p_0]}(L)$ the set whose elements are geodesics
in $D({\Tbar}, \iota)$, of a fixed length $L$ originating
at a fixed point $[1, p_0]$ in ${\Tbar} \subset D({\Tbar}, \iota)$ terminating at a point
in $D({\Tbar}, \iota)$.  We can thus conclude: while ${\cal A}_{p_0}(L)$ is not sequentially
compact, any sequence in ${\cal A}_{p_0}(L)$ converges to an element in
$\tilde{{\cal A}}_{[1,p_0]}(L)$ in the sense that the limiting path is realized in
$D({\Tbar}, \iota)$ as a prolongation of the segment $\overline{p_0 p_1}$ in $\Tbar$, which
could not be prolonged in ${\Tbar}$ as the space is not geodesically complete.
This is also observed by noting that the surjective map ${\rm exp}_{[1,p_0]}^{-1} : D({\Tbar}, \iota)
\rightarrow C_{[1,p_0]} D({\Tbar}, \iota) (\equiv {\bf R}^{6g-6-|\sigma(p_0)|})$ restricted
to ${\Tbar} \subset D({\Tbar}, \iota)$ is no longer surjective.  The difference set
$C_{[1,p_0]} D({\Tbar}, \iota) \backslash {\rm exp}_{[1,p_0]}^{-1} [{\Tbar}]$ consists of
the image of geodesic rays $\theta: {\bf R}_{\geq 0} \rightarrow D({\Tbar}, \iota)$
originating at $p_0$ minus the maximal prolongations
$\overline{p_0 p_1} $ within $\Tbar$ for $\theta$'s.
Recall that $\Tbar$ is a convex subset of $D({\Tbar}, \iota)$.

One should note here that this does not imply that
the development $D({\Tbar}, \iota)$ is locally compact.  For example, consider
a geodesic segment $\overline{p_0 p_1}$ in $\Tbar$ with $p_0 \in {\T}_{\sigma}$ with
$|\sigma| = 1$, $\pi_{\sigma} (p_1) = p_0$ and the distance $d(p_0, p_1)=L$
sufficiently small.  Then the sequence
$\{ \gamma_\sigma^n (\overline{p_0 p_1}) \}$ does not converge in ${\cal A}_{p_0}(L)$, but
it converges to $\overline{[1, p_0][s_\sigma, p_1]}$ in
$\tilde{{\cal A}}_{[1,p_0]}(L)$. Now consider the set
$(\cup_{n=1}^\infty [1, \gamma_\sigma^n p_1]) \cup [s_\sigma, p_1]$ inside the ball
$B_{2L}([1, p_0])$ in $D({\Tbar}, \iota)$.  This set is not compact
as one can choose a covering by sufficiently small balls around each
point of the infinite set $\{[1, \gamma_\sigma^n p_1]\}_{n \in {\bf Z}} \cup [s_\sigma, p_1]$ so that they
are mutually disjoint, for which there is no finite subcover.
Another remark is that the union of the geodesic segments
$\overline{[1,p_1][1, p_0]} \cup \overline{[1, p_0][s_\sigma, p_1]}$
is itself length minimizing.  As a matter of fact, for all $n \in {\bf Z}$,
$\overline{[1,p_1][1, p_0]} \cup \overline{[1, p_0][s(\sigma), \gamma_\sigma^n p_1]}$
is length minimizing.  What this demonstrates is that the exponential map
${\rm exp}_{[1,p_1]}$ is far from being
single-valued beyond the point $[1, p_0]$.

\section{Finite Rank Property}

Given a bounded closed convex set $F$ in a metric space $X$, we define circumradius of $R(F)$ of the set $F$ to be the smallest number
so that there exists a geodesic ball $B_R(O)$ of radius $R$ centered at $O_F$ {\it containing} $F$. The point $O$ is called a circumcenter of $F$.
It is known~\cite{BH} that in a ${\rm CAT}(0)$ space $X$, each bounded closed convex set $F$ is contained in a unique geodesic ball $B_R(O)$ .
For the sake of completeness, we present here a proof of this different from the one in~\cite{BH}.

\begin{theorem}[Existence/Uniqueness of Circumcenter]
For a bounded closed convex set $F$ of circumradius $R$ of an NPC space $X$, there exists a unique circumcenter $O$ in $F$.
\end{theorem}

\begin{proof}
We first define a function on $X$ be
\[
f(x) = \sup_{q \in F} d(x, q).
\]
It is the smallest real number such that the set $F$ is contained in the closed ball of radius $f(x)$ centered at $x$.
Note here that $f(x)$ is a convex functional, being a supremum of bounded convex functionals. Then a circumcenter
would be a point $O$ where $f(x)$ is minimized if it exists.  First note that for any point $x$ in $F$, $f(x)$
is less than or equal to the diameter of the set $F$.  Secondly for any point $x$ outside $F$, we have $f(x) \geq  f(\pi (x))$ where $\pi$
is the nearest projection map of $X$ onto $F$.  Consider the sub-level sets
\[
S_{r} = \{ x \in X: f(x) \leq r \}.
\]
Also define $R(F)$ to be
\[
R(F) = \inf_{x \in X} f(x).
\]
Note that we may take $\inf_{x \in F} f(x)$ to be the value of $R(F)$ instead.
Now we consider a decreasing sequence $\{ R_i \}$ converging to $R$.  Then the sets $\{ S_{R_i} \}$ form a nested
sequence of convex sets, each of which is contained in $F$ for sufficiently large $i$.  Using a fact that a nested decreasing sequence of nonempty closed bounded convex sets in
an ${\rm CAT}(0)$ space has a non-empty intersection (see Proposition 1.2 in~\cite{KS1}) we show here that the diameter
of the nonempty intersection set is zero.  This in turn says that the intersection set consists of a single point,
proving the uniqueness of the circumcenter of the set $F$.

To show that the diameter of the intersection set $\cap_{i \in {\bf N}} S_{R_i}$ is indeed zero,
we chose $P_0$ and $P_1$ to be a pair of points in the convex set $S_{R_i}$ with $P_{1/2}$ their mid-point
and $q$ to be a point in $F$.  We choose a sequence of points $\{q_j^i\}$ in $F$ such that
$\lim_{j \rightarrow \infty} d(P_{1/2}, q_j^i) = R_i$
Now recall the triangle comparison
inequality above with $\lambda = 1/2$:
\[
d(P_{1/2}, q_j^i)^2 \leq \frac{1}{2} d(P_0, q_j^i)^2 + \frac{1}{2} d(P_1, q_j^i)^2 - \frac{1}{4} d(P_0, P_1)^2.
\]
which implies
\[
\frac{1}{4} d(P_0, P_1)^2 \leq \frac{1}{2} d(P_0, q_j^i)^2 + \frac{1}{2} d(P_1, q_j^i)^2 - d(P_{1/2}, q_j^i)^2
\]
Now note that as $i$ and $j$ increase, the right hand side of the inequality above goes to zero, for the distance
$d(P_0, q_j^i)$ and $d(P_1, q_j^i)$ have upper bounds depending on $i$, which are converging down to $R$,
while the distance $d(P_{1/2}, q_j^i)$ converges to $R$ as $i$ and $j$ go to infinity.  By taking the limits,
we obtain that $d(P_0,P_1)=0$ and that the diameter of the set $\cap_i S_{R_i}$ is zero.
\end{proof}

Having established the existence and uniqueness of the circumcenter $O_F$, define {\it circumset} $CS(F)$ of a bounded closed
convex set $F \subset X$ to be $X \cap \partial B_{R(F)} (O_F)$.  We take the closed convex hull ${\rm cvx}(CS (F))$.  This set is
characterized as the smallest closed convex set which has $CS(F)$ as its circumset.  Naturally the diameter $D({\rm cvx}(CS (F)))$ is
less than or equal to $D(F)$, while the circumradius and the circumcenter remain the same.
As the inequality that defines the FR property of a space $X$ is finding a
suitable positive lower bound  for the ratio $D(F)/R(F)$, it suffices to show the inequality for the ratio $D({\rm cvx}(CS (F)))/R(F)$
with the same lower bound.

We now define Alexandrov angle between two geodesic paths
$\sigma$ and $\sigma'$ in a metric space.
\\

\begin{definition}
Let $X$ be an ${\rm CAT}(0)$ space and $\sigma_0$, $\sigma_1$ be constant speed geodesics with
their initial point $p$.
Then Alexandrov angle between the two geodesics is defines by
\[
\angle_p (\sigma_0, \sigma_1) = \lim_{t \rightarrow 0} \arccos \frac{d(p, \sigma_0(t))^2 +
d(p, \sigma_1(t))^2 - d(\sigma_0(t), \sigma_1(t))^2}{2 d(p, \sigma_0(t))d(p, \sigma_1(t))}
\]
\end{definition}

We quote a lemma from~\cite{BH}, which relates the Alexandrov angle with a first variation of the distance function.
\\

\noindent {\bf Proposition}~\cite{BH}(II.3.6) {\it Given a geodesic $\sigma(t)$ with $\sigma(0) = O$ and $y \neq O$
we have}
\[
\cos \angle_O (\sigma(t), y) = \lim_{t \rightarrow 0} \frac{d(\sigma(0), y) - d(\sigma(t), y)}{t}.
\]

Note that this formula implies that $y$ can be replaced by any point $y' \neq O$ on the geodesic segment $\overline{Oy}$
to determine the same angle.  Here the existence of the one-sided derivative is guaranteed by the convexity of the distance
function measured from $p$ and the Lipschitz constant bound, which is equal to one.

In order to characterize the circumcenter from the geometry of the circumset,
we present the following.
\\

\begin{theorem}\label{CC}
Let $\sigma(t)$ be a geodesic originating at a point $O$.  If the point $O$ is the circumcenter of the set $F$
we have the following inequality.
\[
\sup_{q \in CS(F)} \angle_O (\sigma(t), q) \geq \pi/2.
\]
Conversely, if a closed convex hull  ${\rm cvx}(A)$ of a set $A \subset \partial B_R(p)$ satisfies an inequality
\[
\sup_{q \in A} \angle_p (\sigma(t), q) \geq \pi/2
\]
for any geodesic $\sigma$ originating at the point $p$,
then the point $p$ is the circumcenter of the convex hull ${\rm cvx}(A)$.
\end{theorem}

\begin{proof}
We first prove the first half of the statement of the theorem.
Let $\sigma(t)$
be a geodesic originating from the circumcenter $O$ of the set $F$ and hence of the
set $ {\rm cvx}(CS(F))$.
Here we claim that the value of $f(x)=\sup_{q \in
{\rm cvx}(CS(F))} d(x, q)$ equals to the value of
\[
\tilde{f}(x) = \sup_{q \in CS(F)} d(x, q).
\]
Recall the convex hull of a set $CS$ can be constructed inductively.  Namely let $\Omega_0$
be the set $CS$.  Then $\Omega_i$ to be the points expressed as $\lambda x + (1 - \lambda) y$
where $\lambda \in [0,1]$ and $x, y \in \Omega_{i-1}$ for $i= 1,2,\cdots$. Then
the closed convex hull of the circumset is the closure of $\lim_{i \rightarrow \infty}
\Omega_i$. Note that the value of the convex function $f$ restricted to a line segment
$L \cap {\rm cvx }(CS)$ is maximized on an end point of the segment.  As the points of
$\Omega_i \backslash
\Omega_{i-1}$ apprear as the interior points of line segments whose endpoints belong to
$\Omega_{i-1}$, we have
\[
\sup_{x \in\Omega_0} f(x) \geq \sup_{x \in \Omega_i} f(x)  \,\,\, \mbox{ for $i =1,2,\cdots$}
\]
and the statement of the claim follows.

Note that  by definition $\tilde{f}$
is a convex function with its minimum uniquely achieved at $O$, and hence that
there exists non-negative
one-sided limit
\[
\lim_{t \searrow 0} \frac{\tilde{f}(\sigma(t))-\tilde{f}(\sigma(0))}{t} \geq 0
\]
due to the convexity, the Lipschitz constant bound of $d(x, q)$ and the fact
$\tilde{f}(\sigma(t)) \geq \tilde{f}(\sigma(0))=R(F)$.

The hypothesis together with the first variation representation formula of the
Alexandrov angle above,  the inequality we need to show is equivalent to
\[
\sup_{q \in CS} \lim_{t \searrow 0} \frac{d(\sigma(t), q)-d(\sigma(0), q)}{t} \geq 0.
\]
Hence it suffices to justify exchanging the order of taking the supremum and taking
the limit.

Let $\{ q_i \}$ be a sequence in $CS(F)$ so that
\[
\lim_{i \rightarrow \infty} d(\sigma(t), q_i) = \sup_{q \in CS(F)} d(\sigma(t), q) \,\,\,
(= \tilde{f}(\sigma(t)) \,\, )
\]
Then the general theory of convex function~\cite{Ro} would give us the uniform convergence of
the first derivatives of the sequence of the functions to the first derivative
of the limiting function $\tilde{f}(x)$.  Consequently we obtain
\[
\lim_{i \rightarrow \infty} \lim_{t \searrow 0}
\frac{d(\sigma(t), q_i)-d(\sigma(0), q_i)}{t} =
\lim_{t \searrow 0} \frac{\lim_{i \rightarrow \infty} d(\sigma(t), q_i)
- \lim_{i \rightarrow \infty} d(\sigma(0), q_i)}{t}
,
\]
where the right hand side of the equality is, by the argument above, non-negative. Thus
the statement is proven.

As for the second half of the theorem, the hypothesis implies that for all geodesics
$\sigma(t)$ originating at $p$, the first derivative of the function $\sup_{q \in A}
d(\sigma(t), q)$ is non-negative at $t = 0$.  Here we are again exchanging the order of
taking supremum and taking a derivative, which is justified by the convexity of
the functions under considerations.  As the function $\sup_{q \in A} d(\sigma(t), q)$
is convex, it follows that it is monotone increasing for $t \geq 0$. By substituting
the set $A$ by its convex hull ${\rm cvx}(A)$, we have that $\sup_{q \in {\rm cvx}(A)}
d(\sigma(t), q)$ is a monotone increasing convex function. This says that $p$ is the
circumcenter of the set ${\rm cvx}(A)$.

\end{proof}

We now present a proof that the development $D({\Tbar}, \iota)$ and the
\weil completion $\Tbar$ both satisfy  the finite rank property.

\begin{theorem}\label{DFR}
The development $D({\Tbar}, \iota)$ is an FR space.
\end{theorem}

A simple observation is that a complete convex subset of an FR space is again FR.
As the \weil completion $\Tbar$ is a complete convex subset of the development
$D({\Tbar}, \iota)$, we have the statement of the main theorem

\begin{corollary}
The \weil completion $\Tbar$ of a \teich space $\T$ is FR.
\end{corollary}

\begin{proof}(of {\bf Theorem~\ref{DFR}})
We start with a bounded convex set $F$ in the development $D({\Tbar}, \iota)$.
Let $R$ be the circumradius of $F$, $O_F$ be its unique circumcenter.  Recall that the
circumset $CS = CS(F)$ of $F$ is defined to be  $F \cap \partial B_R (O_F)$. Without loss of
generality we assume that $F = {\rm cvx}(CS(F))$, the smallest convex set among all
the convex sets with its circumset $CS(F)$.   Denote by $\overline{CS}$  the image set
$\exp_{O_F}^{-1} [CS] \subset C_{O_F} {D({\Tbar}, \iota)}$.  Lastly denote by $\overline{F}$ the convex hull of
$\overline{CS}$.  Note here that a priori we do not know whether the origin
$\exp_{O_F}^{-1} \{O_F\} $ of the tangent cone $C_{O_F} {\Tbar}$ is contained in $\overline{F}$.

\begin{proposition}\label{OCC}
The origin of the tangent cone $C_{O_F} D({\Tbar}, \iota)$ is the circumcenter of the set $\overline{F}$.
\end{proposition}

In particular, the statement of this proposition implies that the origin is contained in the
convex set $\overline{F}$.

\begin{proof}
Recall that the gluing construction in Section~{\bf 4}, which resulted in modeling the
local picture of the development $D({\Tbar}, \iota)$  so that the tangent cone
$C_p D({\Tbar}, \iota)$ for each $p$ in a strata indexed by an element $(gW_\sigma, \sigma)$
in $D(C({\cal S}), \iota)$ is isometric to a finite
dimensional Euclidean space ${\bf R}^k$, with
$k= \dim {\T}_{\sigma} + |\sigma| = {(6g-6)-2|\sigma|}+|\sigma|$.

The space of directions $S_{O_F} D({\Tbar}, \iota)$ is a $k-1$ dimensional unit sphere
in ${\bf R}^k$.  Let $\sigma(t)$ be a
arc-length parameterized geodesic starting at $O$.  Then the fact that the point $O_F$ is
the circumcenter of the set $F$ implies that
\[
\sup_{q \in CS} \angle_{O_F} (\sigma(t), q) \geq \pi/2.
\]
as the first half of {\bf Theorem~\ref{CC}} guarantees.  Now consider the quotient map
\[
{\rm exp}_{O_F}^{-1}: D({\Tbar}, \iota) \rightarrow C_{O_F} D({\Tbar}, \iota)
\]
where two points $\sigma(t) \sim \sigma'(t)$ of distance $t$ away from $O_F$ are equivalent
when $\angle_{O_F} (\sigma(t), \sigma'(t)) = 0$.  The geodesic extension property of
$D({\Tbar}, \iota)$ guarantees that this map is surjective onto $C_{O_F} D({\Tbar}, \iota)$ isometric
to ${\bf R}^{6g-6-|\sigma|}$.

As the Alexandrov angle is a function defined on rays emanating from the origin of the tangent cone,
the inequality above is rewritten as
\[
\sup_{[q] \in \overline{CS}} \angle_{O} ([\sigma](t)), [q]) \geq \pi/2
\]
where $[\sigma](t)$ is the ray ${\rm exp}_{O_F}^{-1} \sigma(t)$.

As the set $\overline{CS}$ is a subset of $\partial B_{R(F)}(O)$, we in turn apply the second half
of the statement of {\bf Theorem~\ref{CC}}, to conclude that the origin is the circumcenter
of the convex full ${\rm cvx}(\overline{CS}) \subset {\bf R}^k$.

\end{proof}

We quote the following classical result of Caratheodory.
\\

\noindent {\bf Theorem}(Caratheodory~\cite{Ro}) {\it Given a bounded set $S$ in ${\bf R}^k$ and its convex hull $C$, a point $x$ lies in $C$ if and only if
$x$ can be expressed as a convex combination of $k+1$ (not necessarily distinct) vectors $v_i = \overline{xp_i}$, $(1 \leq i \leq k+1)$
with $p_i$ in $S$. }\\

\noindent Apply this result to the situation where the set $S$ is set to be $\overline{CS}$ and $x$ be the origin $O$ in
$C_O D({\Tbar}, \iota) \subset {\bf R}^k$, which is by the previous proposition contained in the convex hull $\overline{F}$ of $\overline{CS}$.  Recall $CS$ is a subset of
$\partial B_{R(F)}(O)$.  Hence $v_1,\cdots,v_{k+1}$ are vectors with $\| v_i \| = R$
for each $i$.  Furthermore, there exists $\lambda_i \in [0,1]$ such that
\[
\sum_{i=1}^{k+1} \lambda_i =1 \mbox{ and  } \sum_{i=1}^{k+1} \lambda_i v_i = O.
\]

Now we claim that for any distinct pair $i$ and $j$, we have
\[
\| v_i - v_j \| \leq D(F)
\]
where $D(F)$ is the diameter of $F$ in $\Tbar$.  To see this, let $p_i$ be in the set $(\exp_{O_F}^{-1})^{-1}\{v_i\} \cap CS$.
Then we first note by definition of diameter
\[
d(p_i, p_j) \leq D(F).
\]
In addition we have
\[
\| v_i - v_j \| \leq d(p_i, p_j)
\]
since the distance between two diverging rays in an NPC space grows super-linearly, in contrast to the linear growth
in Euclidean spaces( see~\cite{BH} (II 3.1), for example.)  Hence combining the two preceding inequalities
results in proving the claim.

Finally to prove the statement of the theorem,
we will show that
\[
R \leq \sqrt{\frac{k}{k+1}} \frac{D}{\sqrt{2}}
\]
where $k$ is the minimal dimension of Euclidean space such that
$C_{O_F} {\Tbar} \subset {\bf R}^k$. Namely there cannot be an isometric embedding of ${\bf R}^l$
into the space $D({\Tbar}, \iota)$ for $l > k$:  if there were a flat of dimension $l > k$
containing the point $O_F$, then its image by ${\rm exp}_{O_F}^{-1}$ in the tangent cone
$C_{O_F} D({\Tbar}, \iota)$ would be isometric to ${\bf R}^l$.

The inequality can be written as
\[
R \leq (1-\varepsilon_0) \frac{D}{\sqrt{2}}
\]
with $\varepsilon_0 = 1 - \sqrt{\frac{k}{k+1}}$, which says that the space is FR ({\bf
Definition~\ref{FR}}), and the
maximal dimension of flats inside the space is $k$.

The inequality is obtained as follows, similar to an argument found in~\cite{KS2};
\begin{eqnarray*}
2 R^2 & = & 2 R^2 \Big( \sum_{i=1}^{k+1} \lambda_i \Big)^2 =  2 R^2 \Big( \sum_{i \neq j} \lambda_i \lambda_j + \sum_k \lambda_k^2 \Big) \\
& = & 2 R^2 \sum_{i \neq j} \lambda_i \lambda_j + 2 \sum_k \lambda_k^2 \| v_k \|^2 \\
& = & 2 R^2 \sum_{i \neq j} \lambda_i \lambda_j - 2 \Big( \| \sum_k \lambda_k v_k \|^2 - \sum_i \lambda_k^2 \| v_k \|^2 \Big) \\
& = & \sum_{i \neq j}  \lambda_i \lambda_j (\|v_i\|^2 + \|v_j\|^2) - 2 \sum_{k \neq l} \lambda_k v_k \cdot \lambda_l v_l \\
& = & \sum_{i \neq j} \lambda_i \lambda_j \| v_i - v_j\|^2 \\
& \leq & \sum_{i \neq j} \lambda_i \lambda_j D^2  =  \Big( (\sum_i \lambda_i)^2 - \sum_j \lambda_i^2 \Big) D^2 \\
& \leq & \Big( 1 - \frac{1}{k+1} \Big) D^2
\end{eqnarray*}
This completes the proof of the finite rank theorem.

\end{proof}

Therefore quoting the result of Korevaar-Schoen~\cite{KS2}, we have the following existence theorem for equivariant harmonic maps into $\Tbar$.

\begin{theorem} Let $\Gamma$ be the fundamental group of a compact Riemannian manifold $M$, and let $\rho$ be an isometric action of
$\Gamma$ on $\Tbar$.  Either there exists an equivalence class of rays, fixed by the
$\rho$-action, or there exists a $\rho$-equivariant harmonic
map $u:\tilde{M} \rightarrow \Tbar$, where $\tilde{M}$ is the universal covering space of $M$.
\end{theorem}

\section{Isometry Group of The Development}

For a Coxeter complex, there is a natural construction of an isometry group of the
Coxeter complex which contains the original Coxeter group as a normal subgroup.
In the current context, the fundamental domain is the \weil completion $\Tbar$ of the \teich
space $\T$, on which the extended  mapping class group $\widehat{{\rm Mod}_{\Sigma}}$ acts isometrically.
The extended mapping class group $\widehat{{\rm Map}_{\Sigma}}$ is known (\cite{Iv},\cite{Ko},\cite{Lu}) to be the full automorphism group of the complex of curves $C({\cal S})$.  Using this fact, it has been
shown~\cite{MW} that the extended mapping class group is indeed the full isometry group
of the \weil completed \teich space.  Note that each element
$\gamma$ of the extended mapping class group $\widehat{{\rm Map}_\Sigma}$ preserves the Coxeter matrix, namely
\[
m_{\gamma(s) \gamma(t)} = m_{st}
\]
As the Coxeter group $W$ is generated by ${\cal S}$, and the group $W$ is
completely determined by the Coxeter matrix $[m_{st}]_{s,t \in {\cal S}}$,
it follows that  each element $\gamma$ in $\widehat{{\rm Mod}_{\Sigma}}$ induces an
automorphism of $W$. Such an automorphism of $W$ is called {\it diagram automorphism}
\cite{Da}.

The formalism laid out in M.Davis' book gives us a natural action (Proposition 9.1.7~\cite{Da}) of the
semi-direct product $G:=W \rtimes \widehat{{\rm Mod}_\Sigma}$ on the development
$D({\Tbar}, \iota)$ as follows: given $u=(g, \gamma) \in G$ and
$[g', y] \in D({\Tbar}, \iota)$,
\[
u \cdot [g', y] := [g \gamma (g'), \gamma y]
\]
where $\gamma(g')$ is the image of $g'$ by the automorphism of $W$ induced by
$\gamma:C({\cal S}) \rightarrow C({\cal S})$.

Clearly the action $G \hookrightarrow D({\Tbar}, \iota)$ is isometric. Thus $G$ is a
subgroup of the isometry group ${\rm Isom}(D({\Tbar}, \iota))$.  It remains an open question
that if this group is indeed the full isometry group.

We note also the following.

\begin{theorem}
The action of the Coxeter group $W \subset G $
on $D({\Tbar}, \iota)$ is properly discontinuous.
\end{theorem}

A group $\Gamma$ acts on a Hausdorff
space $Y$ properly discontinuously if i) $Y/\Gamma$ is Hausdorff, ii) for each $y \in Y$, the isotropy
subgroup $\Gamma_y =\{ g \in \Gamma| gy=y\}$ is finite, and iii) Each $y \in Y$ has a $\Gamma_y$-stable
neighborhood $U_y$ such that $gU_y \cap U_y = \emptyset$ for all $g \in \Gamma \backslash \Gamma_y$.

A statement (Lemma~5.1.7~\cite{Da}) tells us that in our context where
 $\Gamma = W_{\cal S}$ and $Y = D({\Tbar}, \iota)$, we only need to check that for any subset $\hat{\cal S} \subset {\cal S}$ such that
the subgroup $W_{\hat{\cal S}}$ generated by $\hat{\cal S}$ is of infinite order, the fixed point set of $W_{\hat{\cal S}}$ is empty.
This condition is clearly satisfied here, because $\hat{\cal S}$ is not an element of the
 complex of curves $C({\cal S})$, and thus
there is no stratum indexed by $\hat{\cal S}$.

Lastly we can quote the
theorem by Korevaar-Schoen~\cite{KS2} again in this context to obtain the existence
theorem of harmonic map

\begin{theorem} Let $\Gamma$ be the fundamental group of a compact Riemannian manifold $M$,
and let $\rho: \Gamma \rightarrow G$ be a representation which defines an isometric action of
$\Gamma$ on $D({\Tbar}, \iota)$.  Then either there exists an equivalence
class of rays, fixed by the $\rho$-action, or there exists a $\rho$-equivariant harmonic map
$u:\tilde{M} \rightarrow D({\Tbar}, \iota)$, where $\tilde{M}$ is the universal covering
space of $M$.
\end{theorem}

Once again this dichotomy regarding the existence of harmonic map into the development $D({\Tbar}, \iota)$
is the exact analogue of the much-investigated situation for harmonic map into non-compact symmetric
spaces $G/K$ for semi-simple Lie group $G$.

\end{document}